\documentclass[11pt,a4paper]{article}

\usepackage[rgb]{xcolor}
\usepackage{hyperref}
\usepackage[utf8]{inputenc}
\usepackage[english]{babel}
\usepackage[babel]{csquotes}
\usepackage{fullpage}
\usepackage{amsfonts, amsthm, amsmath, amssymb}
\usepackage{bbm}       
\usepackage[version=3]{mhchem}    
\usepackage{tikz}      
\usetikzlibrary{shapes,automata,positioning,arrows}
\usepackage{mathtools}  

\numberwithin{equation}{section}

\newtheorem{theorem}{Theorem}[section]
\newtheorem{corollary}[theorem]{Corollary}
\newtheorem{lemma}[theorem]{Lemma}
\newtheorem{proposition}[theorem]{Proposition}

\theoremstyle{definition}
\newtheorem{definition}{Definition}[section]

\newtheorem{example}{Example}[section]

\def\SS{\mathcal{S}}
\newcommand{\C}{\mathcal{C}}
\newcommand{\R}{\mathcal{R}}
\newcommand{\G}{\mathcal{G}}

\newcommand{\NN}{\mathbb{N}}
\newcommand{\RR}{\mathbb{R}}
\newcommand{\ZZ}{\mathbb{Z}}

\def\wt{\widetilde}

\definecolor{mygrey}{gray}{0.75}

\newcommand{\been}{\begin{enumerate}}
\newcommand{\enen}{\end{enumerate}}
\newcommand{\beit}{\begin{itemize}}
\newcommand{\enit}{\end{itemize}}

\def\one{\mathbbm{1}}

\DeclareMathOperator{\SSann}{span}
\DeclareMathOperator{\kernel}{Ker}
\DeclareMathOperator{\supp}{supp}

\providecommand{\abs}[1]{\lvert#1\rvert}
\providecommand{\norm}[1]{\lVert#1\rVert}

\title{Transition graph decomposition for complex balanced reaction networks with non-mass-action kinetics}
\author{Daniele Cappelletti\footnotemark[1] \and Badal Joshi\footnotemark[2]}

\begin{document}

 \footnotetext[1]{Department of Mathematics, Politecnico di Torino, Italy}
  \footnotetext[2]{Department of Mathematics, California State University San Marcos}

 \tikzset{every node/.style={auto}}
 \tikzset{every state/.style={rectangle, minimum size=0pt, draw=none, font=\normalsize}}
  \tikzset{bend angle=7}

 \maketitle

\begin{abstract}
Reaction networks are widely used models to describe biochemical processes. Stochastic fluctuations in the counts of biological macromolecules have amplified consequences due to their small population sizes. This makes it necessary to favor stochastic, discrete population, continuous time models. The stationary distributions provide snapshots of the model behavior at the stationary regime, and as such finding their expression in terms of the model parameters is of great interest. The aim of the present paper is to describe when the stationary distributions of the original model, whose state space is potentially infinite, coincide exactly with the stationary distributions of the process truncated to finite subsets of states, up to a normalizing constant. The finite subsets of states we identify are called {\em copies} and are inspired by the modular topology of reaction network models. With such a choice we prove a novel graphical characterization of the concept of complex balancing for stochastic models of reaction networks. The results of the paper hold for the commonly used mass-action kinetics but are not restricted to it, and are in fact stated for more general setting.
\end{abstract}

\section{Introduction}\label{sec:introduction}

Reaction networks are mathematical models for studying the evolution of biochemical systems, widely used in applications. Formally, they are made of a set of species whose populations evolve in time guided by a set of reactions such as the ones in the example below:
 \begin{equation}\label{eq:reactions}
 A+B \cee{<=>[\kappa_1][\kappa_2]} 2C ~,~  A \cee{<=>[\kappa_3][\kappa_4]} B. 
  \end{equation}
Each species has a discrete number of molecules whose counts change as the molecules react or are produced from reactions. 
In many circumstances, the counts are sufficiently high so that the granularity may be ignored and it suffices to model species concentrations rather than counts. 
However, several biological macromolecules such as proteins and enzymes are present only in small numbers ($\sim 100$) and it becomes necessary to both track individual numbers and be mindful of stochastic fluctuations in those numbers.
Due to such considerations, $\ZZ^n$ is the natural state space for the stochastic model (a continuous-time Markov chain) that we consider here. 

The stationary distribution of a stochastic system describes the state of an isolated system over a long period of time. 
So in a sense it is an attracting steady state of the dynamic model. 
Having an analytic expression of the stationary distribution in terms of the model parameters is useful for various reasons: first of all, the expected level and the fluctuations of different proteins in stationary biological systems is of great interest. Controlling the stationary regime is also attracting an increasing interest, given the novel technological possibilities of changing the cell DNA and, hence, the rate and the form of the biological reactions occurring inside. Notable examples of how the emerging field of synthetic biology focuses on the control of stationary regimes are given in \cite{integral, plesa2021quasi, kim2020absolutely}. Finally, the approximation of complex biological systems with reaction rates spanning over different orders of magnitudes can be performed if the stationary distribution of the faster subsystem is known \cite{kang:separation}. The study of multiscale biological system is of crucial importance and is recognized as one of the key ingredients to unlock the function of genetic and cellular compositions in the Perspectives in Mathematical Biology individuated by the European Society for Mathematical and Theoretical Biology (ESMTB) \cite{challenges, multiscale}.

The aim of the present paper is twofold: first, typical computational techniques to calculate an approximation of the stationary distribution of interest consists in truncating the (potentially infinite) state space to a smaller finite set, and calculate the stationary distribution of the process restricted to it \cite{fsp, trunc}. The approximation error is often hard to calculate. A natural question is whether the same stationary distribution of the full model is also a stationary distribution of the restricted one; if that is the case then the approximation error only concerns the normalizing constant, but not the ratio of the distribution in different states. In the present paper we consider restrictions of the state space to the images of {\em copies} (or union thereof). A copy of a reaction network is an embedding of the reaction graph into the state space. Figure~\ref{fig} illustrates how the state space of
\begin{equation}\label{eq:ex}
  \begin{tikzpicture}[baseline={(current bounding box.center)}]
   \node[state] (0)   at (0,0)    {$0$};
   \node[state] (A+B)   at (1,1)    {$A+B$};
   \node[state] (A)   at (2,0) {$A$};
   \path[->]
    (0) edge node {} (A+B)
    (A+B) edge node {} (A)
    (A) edge node {} (0);
  \end{tikzpicture}
\end{equation}
can be decomposed into copies. 

\begin{figure}[h]
\label{fig}
\begin{center}
   \begin{tikzpicture}
   \draw[->,>=stealth,shorten <=0.3cm,shorten >=0.3cm, line width=0.5mm, draw=mygrey] (-1,0) to (7.5,0);
   \node[text=gray] at (-.5,3.2) {B};
    \node[text=gray] at (7.7,-.1) {A};
   \clip (-0.5,-0.5) rectangle (6.8cm,3.8cm); 
 \draw[->,>=stealth,shorten <=0.3cm,shorten >=0.3cm, line width=0.5mm, draw=mygrey] (0,-1) to (0,4);
 \foreach \x in {0,1,...,7}{                           
    \foreach \y in {0,1,...,5}{                       
    \pgfmathsetmacro{\huenum}{mod(6*\x+20*\y,48)/48};
    \definecolor{mycolor}{hsb}{\huenum,1,1};
    \node[circle,fill=black, inner sep=2pt, anchor=center] at (\x,\y) {}; 
    \draw[->,>=stealth,thick,shorten <=0.3cm,shorten >=0.3cm, draw=mycolor] (\x,\y) to node {} (\x+1, \y+1);
    \draw[->,>=stealth,thick,shorten <=0.3cm,shorten >=0.3cm, draw=mycolor] (\x+1,\y+1) to node {} (\x+1, \y);
    \draw[->,>=stealth,thick,shorten <=0.3cm,shorten >=0.3cm, draw=mycolor] (\x+1,\y) to node {} (\x, \y);
    }
}
\end{tikzpicture}
\end{center}
\caption{The figures depicts the transition graph of the continuous time Markov chain associated with the reaction network \eqref{eq:ex}. Different colors indicate different {\em copies}, which are different embedding of the graph \eqref{eq:ex} in the state space $\ZZ^2_{\geq0}$. }
\end{figure}

{\em We prove that the process restricted to arbitrary unions of copies has the same stationary measure of the original one only in the case of complex balancing}. In the case of complex balanced systems, an exact formula for the stationary distributions is already known \cite{anderson:product-form}, however the general line of thought can be pursued for other types of restrictions different from copies. In fact, in future work we will use this idea to efficiently calculate novel stationary distributions. 

Secondly, the paper focuses on the connection between stationary distributions and graphical symmetries. This perspective is based on links between steady states and graphical features of the model which have been fruitfully and extensively studied for deterministic models of reaction networks \cite{horn1972general,feinberg1972,sturmfels,gopalkrishnan2014geometric,joshi2013atoms,joshi2022minimal, perez, cappelletti2020hidden, craciun2021multistationarity}. Connections between graphical symmetries of the model and the shape of its stationary distribution have been successfully unlocked in several recent papers \cite{anderson:product-form,cappelletti:complex_balanced,joshi:detailed,CJ2018, hoessly}. In this work, we define a new type of graphical symmetry (balancing of copies, as defined below) and we show that it characterizes the notion of complex balancing under general assumptions. 

While our study is valid for the commonly used mass-action kinetics, the scope of our results is not restricted to such kinetics but is much broader.
Some general results can be stated for {\em any admissible choice of reaction rates}, and others hold for a generalization of mass-action kinetics that we call {\em product form kinetics}, considered for example in \cite{anderson:product-form, AN2018, kelly1979reversibility, whittle1986systems}.

%

\section{Background}
 
 \subsection{Notation} 
 
Let $\RR$, $\RR_{\ge 0}$ and $\RR_{> 0}$ represent the reals, the non-negative reals and the positive reals, respectively. Let $\ZZ$, $\ZZ_{\ge 0}$ and $\ZZ_{> 0}$ represent the integers, the non-negative integers and the positive integers, respectively. 
For $v\in\RR^n$, $\norm{v}_1=\abs{v_1 + \ldots + v_n}$.
 For $v,w\in\RR^n$, $v \le w$ ($v < w$) means that $v_i \le w_i$ ($v_i < w_i$) for all $i \in \{1,\ldots, n\}$. 
For $v,w\in\RR^n$, we define  
 \[
 \one_{\{v \le w\}} = 
 \begin{cases}
 1 \quad , \quad v \le w \\
 0 \quad , \quad \mbox{ otherwise}.
 \end{cases}
 \]
If $v>0$ then $v$ is said to be positive. Similarly, if $v\geq0$ then $v$ is said to be nonnegative. Finally,  $\supp v$ denotes the index set of the non-zero components. For example, if $v=(0,1,1)$ then $\supp v=\{2,3\}$.
 
 If $x\in\RR^n_{\ge 0}$ and $v\in \ZZ^n_{\ge 0}$, we define
 $$x^v=\prod_{i=1}^n x_i^{v_i},\quad \text{and}\quad v!=\prod_{i=1}^n v_i!,$$
 with the conventions that $0!=1$ and $0^0=1$. 

 \subsection{Reaction networks}  
 
 A reaction network is a triple $\G=(\SS,\C,\R)$, where $\SS$ is a set of $n$ species, $\C$ is a set of $m$ complexes, and $\R\subseteq \C\times\C$ is a set of $r$ reactions, such that  $(y,y)\notin\R$ for all  $y\in\C$. The complexes are linear combinations of species over $\NN$, identified as vectors in $\RR^n$.  A reaction $(y,y')\in\R$ is denoted by $y\to y'$, and the vector $y'-y$ is the corresponding \emph{reaction vector}. We require that every species has a nonzero coordinate in at least one complex and that every complex appears in at least one reaction. With this convention, there are no ``redundant'' species or complexes and $\G$ is uniquely determined by $\R$. In \eqref{eq:reactions}, there are $n=3$ species ($A,B,C$), $m=2$ complexes ($A+B, 2B$), and $r=2$ reactions. 
 
 A reaction network $\G=(\SS,\C,\R)$ can be viewed as a graph with node set $\C$ and edge set $\R$ in a natural manner. We will frequently use the viewpoint of a reaction network as a graph in the rest of the paper.

A reaction network $\G$ is \emph{weakly reversible} if every reaction $y\to y'\in\R$ is contained in a closed directed path. Moreover, $\G$ is \emph{reversible} if for any reaction $y\to y'\in\R$, $y'\to y$ is in $\R$. It is clear that each reversible reaction network is also weakly reversible. As an example, the network in \eqref{eq:reactions} is reversible, and therefore weakly reversible. 
 
  The \emph{stoichiometric subspace} of $\G$ is the linear subspace of $\RR^n$ generated by the reaction vectors, namely
 $$S=\SSann(y'-y\vert y\to y'\in\R).$$
 For $v\in\RR^n$, the sets $(v+S)\cap\RR^n_{\ge 0}$ are called the \emph{stoichiometric compatibility classes} of $\G$. 
 
 \subsection{Reaction systems}
 We will consider dynamics of a reaction network with $n$ species both on $\RR^n_{\ge 0}$ and $\ZZ^n_{\ge 0}$. $\RR^n_{\ge 0}$ is the usual underlying state space for deterministic models, while for classic stochastic models the state space is $\ZZ^n_{\ge 0}$. We do not consider stochastic differential equations (ordinary differential equations with a noise term) in this paper, but in passing we mention that this is an instance where a stochastic model has the underlying state space $\RR^n_{\ge 0}$, see for instance \cite{kurtz1976, kurtzstrong, ruth, bibbona:weak}.
  
%

 \subsubsection{Deterministic dynamics}
Let $\G$ be a reaction network. We want to associate each reaction $y \to y' \in \R$ with a {\em rate function} $\lambda_{y \to y'}$, whose domain is $\RR^n_{\ge 0}$. Formally, we define as {\em deterministic kinetics} the following correspondence between reactions and rate functions: 
\[
\Lambda: (y \to y') \mapsto \lambda_{y \to y'}
\]
It is required that for each $y\to y'\in\R$
\begin{equation}\label{eq:cond_det}
 \lambda_{y\to y'}(x)>0\quad\text{only if }\supp y\subseteq \supp x.
\end{equation}
The pair $(\G,\Lambda)$ is called {\em continuous reaction system}. In the deterministic context, the evolution of the species concentrations $z(t) \in \RR_{\ge 0}^n$ is determined by the system of ODEs
 \begin{equation}\label{eq:ODE}
  \frac{d z}{dt}=\sum_{y\to y'\in\R}(y'-y)\lambda_{y\to y'}(z).
 \end{equation}
 Condition \eqref{eq:cond_det} implies that any solution to $\eqref{eq:ODE}$ is non-negative at all times for which it is defined. The solution $z(t)$ is also confined to its {\em stoichiometric compatibility class}:
 \[
 z(t)\in (z(0)+S)\cap\RR^n_{\ge 0}. 
 \]
A state $c \in \RR^n$ is said to be a {\em steady state} of a continuous reaction system $(\G,\Lambda)$ if 
\[
\sum_{y\to y'\in\R}(y'-y)\lambda_{y\to y'}(c) = 0.
\]

\subsubsection{Stochastic dynamics}

Let $\G$ be a reaction network. We now want to associate each reaction $y \to y' \in \R$ with a {\em rate function} $\lambda_{y \to y'}$, whose domain is $\ZZ^n_{\ge 0}$. Similarly to before, we define as {\em stochastic kinetics} the following correspondence between reactions and rate functions: 
\[
\Lambda: (y \to y') \mapsto \lambda_{y \to y'}.
\]
We require
\begin{equation}\label{eq:stoch_cond}
 \lambda_{y\to y'}(x)>0\quad\text{only if }x\geq y.
\end{equation}
The pair $(\G,\Lambda)$ is called {\em stochastic reaction system}. In this setting, a state $x = (x_1 ,\ldots, x_n) \in \ZZ^n_{\ge 0}$ represents the counts $x_i$ of each species $i = 1, \ldots, n$. $X(t)$ represents the state of the system at time $t$, and is considered to be a continuous-time Markov chain with transition rate from state $x$ to state $x'$ given by
\begin{equation*}
 q(x,x')=\sum_{\substack{y\to y'\in\R\\ y'-y=x'-x}} \lambda_{y\to y'}(x).
\end{equation*}
Equivalently, each rate function described the rate of occurrence of the associated reaction, and whenever a reaction $y\to y'$ takes place the process $X(\cdot)$ moves from the current state $x$ to $x+y'-y$. Condition \eqref{eq:stoch_cond} forces the process $X(\cdot)$ to the positive orthant.

We say that $x' \in \ZZ^n_{\ge 0}$ is {\em accessible} from $x \in \ZZ^n_{\ge 0}$, if there is a sequence of states $(x=u_0, u_1, \ldots, u_{n-1}, u_n = x')$ such that for each consecutive pair of states $(u_i, u_{i+1})$, $0 \le i \le n-1$, we have $q(u_i, u_{i+1})>0$.  
A non-empty set $\Gamma\subseteq\ZZ_{\ge 0}^{n}$ is an \emph{closed irreducible set} of $(\G,\Lambda)$ if for all $x\in\Gamma$ and all $u\in\ZZ_{\ge 0}^{n}$, $u$ is accessible from $x$ if and only if $u\in\Gamma$ \cite{norris:markov}. A probability distribution $\pi$ is a stationary distribution for a continuous-time Markov chain if for all $t\in\RR_{\geq0}$ and $x\in\ZZ^n_{\geq0}$
\begin{equation*}
 P(X(t)=x|X(0)\sim \pi)=\pi(x).
\end{equation*} 
By standard theory on Markov chains, the support of $\pi$ (that is, the largest set of states where $\pi$ has a positive value) is a union of closed irreducible sets \cite{norris:markov}. Moreover, by standard theory on Markov chains, if the continuous-time Markov chain is non-explosive (in the sense of \cite{norris:markov}) then a probability distribution $\pi$ is stationary if and only if for all states $x$
\begin{equation}\label{eq:master}
 \pi(x) \sum_{y\to y' \in \R}  \lambda_{y \to y'}(x) = \sum_{y\to y' \in \R} \pi(x+y - y') \lambda_{y \to y'}(x+y - y').
\end{equation}
We further give the definitions below.
\begin{definition} 
Let $(\G,\Lambda)$ be a discrete reaction system and let $\nu$ be a measure on $\ZZ^n_{\ge 0}$. 
\been
\item $\nu$ is said to be a {\em stationary measure} if for all $x \in \ZZ^n_{\ge 0}$
\begin{equation}\label{eq:master2}
 \nu(x) \sum_{y\to y' \in \R}  \lambda_{y \to y'}(x) = \sum_{y\to y' \in \R} \nu(x+y - y') \lambda_{y \to y'}(x+y - y').
\end{equation}
\item $\nu$ is said to be a {\em $\sigma$-finite measure} if $\nu(x) < \infty$ for all $x \in \ZZ^n_{\ge 0}$. $\pi$ is said to be a {\em finite measure} if $\sum_{x\in \ZZ^n_{\ge 0}}\nu(x) < \infty$. 
\enen
\end{definition}

Note that based on the definitions used in this paper, if a stationary measure is a probability distribution is not necessarily a stationary distribution, unless the model is non-explosive. Non-explosiveness of the process can be assessed if it is {\em complex balanced}, as stated in Theorem~\ref{thm:finite_non_explosive} below.

\subsubsection{Mass action kinetics}
An important choice of kinetics is {\em mass-action kinetics}. 
\begin{definition} Consider a reaction network $\G = (\SS,\C,\R)$. 
\been
\item A {\em deterministic mass-action system} is a continuous reaction system $(\G,\Lambda_\kappa^D)$ with 
\[
K_\kappa^D(y\to y')(z)=\kappa_{y\to y'} z^{y}
\]
for some constants $\kappa_{y \to y'}\in \RR_{>0}$, called {\em rate constants}. 
\item A {\em stochastic mass-action system} is a discrete reaction system $(\G,K_\kappa^S)$ with
\begin{equation}\label{eq:sma}
K_\kappa^S(y\to y')(x)=\kappa_{y\to y'}\frac{x!}{(x-y)!}\mathbbm{1}_{\{x\geq y\}}
\end{equation}
for some constants $\kappa_{y \to y'}\in \RR_{>0}$, called {\em rate constants}.  
\enen
\end{definition}  
Many generalizations of mass-action kinetics have been proposed both for stochastic and deterministic models. In the present paper, we will study the following generalization for stochastic models, considered for example in \cite{anderson:product-form, AN2018, kelly1979reversibility, whittle1986systems}. We define the {\em stochastic product form kinetics} as the discrete kinetics
\begin{equation}\label{eq:general}
 K^S_{\kappa, \theta}(y\to y')(x)=\kappa_{y\to y'}\prod_{i=1}^n \prod_{j=0}^{y_i-1}\theta_i(x_i-j),
\end{equation}
for some constants $\kappa_{y \to y'}\in \RR_{>0}$, called {\em rate constants}, and some functions $\theta_i\colon \ZZ\to \RR_{\geq0}$ satisfying $\theta_i(m)=0$ if and only if $m\leq0$. A stochastic reaction system $(\G,K_{\kappa,\theta}^S)$ will be called {\em stochastic product form system}. 

Note that if $\theta_i(m)=m\mathbbm{1}_{\{m\geq0\}}$ for all $i$, then $K^S_{\kappa, \theta}=K^S_\kappa$. Finally, we give the following definition:
\begin{definition}
 A stochastic product form system is called {\em non-saturating} if $\lim_{m\to\infty}\theta_i(m)=\infty$ for all $1\leq i\leq n$.
\end{definition}

\section{Complex balancing}

The definition of complex balancing dates back to \cite{feinberg1972, horn1972general}, where deterministic models were of interest. The definition is the following:
\begin{definition}
 Let $(\G, \Lambda)$ be a continuous reaction system. A state $c\in\RR^n_{\geq0}$ is \emph{complex balanced} if for all complexes $y\in\C$ we have
 \begin{equation*}
  \sum_{y'\in\C\,:\,y\to y'\in\R} \lambda_{y\to y'}(c)=\sum_{y'\in\C\,:\,y'\to y\in\R} \lambda_{y'\to y}(c).
 \end{equation*}
\end{definition}
In words, $c$ is complex balanced if the sum of rates of reactions ``entering'' any complex $y$ is equal to the sum of the rates of reactions ``exiting'' from $y$, calculated in $c$. It is known that every complex balanced state is a steady state, while the opposite is not true in general \cite{horn1972general}. The following are classical results from \cite{horn1972general, feinberg1972}:
\begin{theorem}\label{thm:classic}
 Let $(\G, K^D_\kappa)$ be a deterministic mass-action system. If a positive complex balanced state $c$ exists, then $\G$ is weakly reversible, all positive steady states are complex balanced, and exactly one positive steady state exists within each stoichiometric compatibility class. 
\end{theorem}
Theorem~\ref{thm:classic} implies that the positive steady states of a deterministic mass-action system are either all complex balanced, or none of them is. Moreover, in the former case non-positive steady states are necessarily complex balanced as well as shown in \cite[Theorem 4]{cappelletti:complex_balanced}. Mass action systems are called {\em complex balanced} if at least one positive steady state is complex balanced, which by Theorem~\ref{thm:classic} is equivalent to the existence of at least a positive steady state and to  all of them being complex balanced.

In \cite{anderson:product-form} the first connection between complex balancing and stochastic dynamics is performed, and the following result is proven.
\begin{theorem}\label{thm:productform}
 Let $\G$ be a reaction network. Assume the deterministic mass-action system $(\G, K^D_{\kappa})$ is complex balanced and $c\in\RR^n_{>0}$ is a positive complex balanced state. Then the stochastic reaction system $(\G, K^S_{\kappa, \theta})$ with $K^S_{\kappa, \theta}$ defined as in \eqref{eq:general} and the same rate constants as $(\G, K^D_{\kappa})$ has a stationary measure of the form
 \begin{equation}\label{eq:general_sm}
  \nu(x)=c^x\prod_{i=1}^n\prod_{j=1}^{x_i}\frac{1}{\theta_i(j)}
 \end{equation}
\end{theorem}
If $K^S_{\kappa, \theta}$ reduced to mass-action kinetics, then \eqref{eq:general_sm} becomes proportional to a product-form Poisson distribution, which is a stationary distribution because the process is non-explosive as proven in \cite{ACK:explosion}. In \cite[Theorem 4.1]{AN2018} mild conditions implying finiteness of \eqref{eq:general_sm} and non-explosiveness of the model are derived:
\begin{theorem}\label{thm:finite_non_explosive}
 Let $\G$ be a reaction network. Assume the deterministic mass-action system $(\G, K^D_\kappa)$ is complex balanced and $c\in\RR^n_{>0}$ is a positive complex balanced state. Consider a stochastic non-saturating product form system $(\G, K^S_{\kappa, \theta})$ with the same rate constants as $(\G, K^D_\kappa)$. Then, the associated process $X(\cdot)$ is non-explosive and the measure \eqref{eq:general_sm} is finite.
\end{theorem}
Theorem~\ref{thm:classic}, non-explosiveness of the process and finiteness of \eqref{eq:general_sm} imply that \eqref{eq:general_sm} is proportional to a stationary distribution of $X(\cdot)$. In \cite{cappelletti:complex_balanced} a parallel theory on complex balancing is developed for stochastic models. A complex balanced measure is defined as follows.
\begin{definition}\label{def:sm}
 Let $(\G, \Lambda)$ be a discrete reaction system. A measure $\nu$ with support in $\ZZ^n_{\geq0}$ is \emph{complex balanced} if for all complexes $y\in\C$ and all states $x\in\ZZ^n_{\geq0}$ we have
 \begin{equation}\label{eq:cb_def}
  \sum_{y'\in\C\,:\,y\to y'\in\R} \lambda_{y\to y'}(x)\nu(x)=\sum_{y'\in\C\,:\,y'\to y\in\R} \lambda_{y'\to y}(x+y'-y)\nu(x+y'-y).
 \end{equation}
\end{definition}
In words, a measure is complex balanced if, at any state $x$, the sum of the ingoing fluxes of reactions ``entering'' any complex $y$ is equal to the sum of the outgoing fluxes of reactions ``exiting'' from $y$. It is not difficult to check that a complex balanced measure is necessary stationary, while the converse does not hold in general. Definition~\ref{def:sm} was given for probability distributions in \cite{cappelletti:complex_balanced} and then extended to measures in \cite{CJ2018}. In particular, it is proven in \cite[Proposition 4.13]{CJ2018} that under the hypothesis of mass-action a $\sigma$-finite complex balanced measure is finite, hence normalizable to a probability distribution. 
A strong connection between complex balancing in stochastic and deterministic models exists, as stated below.
\begin{theorem}\label{thm:cbiffcb}
 A stochastic non-saturating product form system $(\G,K^S_{\kappa, \theta})$ has a $\sigma$-finite, positive, complex balanced measure if and only if the deterministic mass-action system $(\G,K^D_\kappa)$ with the same choice of rate constants is complex balanced.
\end{theorem}
Theorem~\ref{thm:cbiffcb} is proven in \cite{cappelletti:complex_balanced} in the case of stochastic mass-action kinetics. The proof for the more general result above is given in Section~\ref{sec:proofs}.

\section{Main results}
 Here we state and prove our main results.
\begin{definition}\label{def:copy}
Let $\G$ be a reaction network. A function $f:\C \to\ZZ_{\ge 0}^n$ is a  {\em copy of $\G$} if for every $y\to y'\in\R$, we have $f(y')-f(y)=y'-y$. We further say that a copy $f$ is {\em active} if every $y\to y'\in\R$ is active at $f(y)$.
\end{definition}
Denote by $f(\C)$ the image of $f$. $\G$ induces a directed graph on $f(\C)$ in a natural manner by associating a directed edge from $f(y)$ to $f(y')$ whenever $y \to y' \in \R$. 
With this understanding, the inclusion copy $\iota: \C \hookrightarrow \ZZ^n_{\ge 0}$ corresponds to the geometrically embedded graph of $\G$ defined in \cite{craciun2015toric}. We state here a lemma that will be useful to prove our main results.
\begin{lemma}\label{lem:inj_copy_exists}
 Let $(\G,\Lambda)$ be a stochastic reaction system. Let $y\to y'$ be a reaction active at a state $x$. Then, there exists an injective copy $f$ of $\G$ with $f(y)=x$. Moreover, if $\Lambda$ is such that $\lambda_{y\to y'}(\tilde{x})>0$ if and only if $\tilde{x}\geq y$, then there exists an injective, active copy $f$ of $\G$ with $f(y)=x$.
\end{lemma}
\begin{proof}
 Let $h=x-y$. Since $y\to y'$ is active at $x$, it follows from \eqref{eq:stoch_cond} that $h\in\ZZ^n_{\geq0}$. Consider the function $f\colon\C\to \ZZ_{\ge 0}^n$ defined by $f(\wt y)=h+\wt y$. We have $f(y)=h+y=x$. Moreover, for all two complexes $\wt y, \wt y'\in\C$ we have $f(\wt y')-f(\wt y)=\wt y'-\wt y$, so $f$ is a copy of $\G$ and is injective. Moreover, if $\lambda_{y\to y'}(\tilde{x})>0$ if and only if $\tilde{x}\geq y$ then $f$ is active. Hence, the proof is concluded.
\end{proof}

\begin{definition}
 Let $(\G,\Lambda)$ be a stochastic reaction system. Let $f$ be a copy of $\G$ and let $\nu$ be a measure on $\ZZ^n_{\ge 0}$. 
 \beit
\item $f$ is {\em active with respect to $(\Lambda, \nu)$} if for all $y \to y' \in \R$
    \begin{equation*}
 \nu(f(y)) \sum_{\substack{\wt y\to \wt y'\in\R:\\f(\wt y) = f(y), f(\wt y') = f(y')}}  \lambda_{\wt y\to \wt y'}(f(y)) > 0
  \end{equation*}
 
\item $f$ is {\em node balanced with respect to $(\Lambda, \nu)$} if for every $x \in f(\C)$
   \begin{equation}\label{eq:nb_def}
 \nu(x) \sum_{y\to y'\in\R:f(y) = x}  \lambda_{y\to y'}(x) =  \sum_{y'\to y\in\R:f(y) = x} \nu(f(y')) \lambda_{y'\to y}(f(y'))
  \end{equation}
  \enit
 \end{definition}

To understand the meaning of node balancing, let $W_f(\cdot)$ be the continuous-time Markov chain whose transition graph is the one induced $f$, and the transition rates are given by
\begin{equation*}
 q(x,x')=\sum_{\substack{y\to y'\in\R:\\f(y)=x, f(y')=x' }}\lambda_{y\to y'}(x).
\end{equation*}
Then, the following holds.
\begin{proposition}\label{prop:copy_bal}
 Let $(\G,\Lambda)$ be a stochastic reaction system. Let $f$ be a copy of $\G$ and let $\nu$ be a $\sigma$-finite measure of $\ZZ^n_{\ge 0}$. Then, there exists a stationary distribution of $W_f(\cdot)$ proportional to the restriction of $\nu$ to $f(\C)$ if and only if $f$ is node balanced with respect to $(\Lambda, \nu)$.
\end{proposition}
\begin{proof}
 The state space of $W_f$ is $f(\C)$, which is finite. Then, there exists a stationary distribution proportional to the restriction of $\nu$ to $f(\C)$ if and only if the latter is a stationary measure, which is equivalent to \eqref{eq:nb_def}. The proof is then concluded.
\end{proof}
Proposition~\ref{prop:copy_bal} implies that if a copy $f$ is node balanced with respect to $(\Lambda, \pi)$, then the restriction of $\pi$ to $f(\C)$ is stationary for $W_f$. If $\pi$ is a stationary distribution of $X(\cdot)$, then $\pi$ is stationary for both $X(\cdot)$ and $W_f(\cdot)$. The same holds for union of copies, as expressed in the following corollary.
\begin{corollary}\label{cor:copy_bal}
Let $(\G,\Lambda)$ be a stochastic reaction system and let $\nu$ be a $\sigma$-finite measure of $\ZZ^n_{\ge 0}$. 
For some $h \ge 1$, let $\Phi =\{f_1,\ldots, f_h\}$ be a set of node balanced copies of $\G$ with respect to $(\Lambda, \nu)$.
Then, the continuous-time Markov chain on state space $\ZZ^n_{\ge 0}$ with transition rates 
\begin{equation}\label{eq:rate}
 q(x,x')= \sum_{i=1}^h \sum_{\substack{y\to y'\in\R:\\f_i(y)=x, f_i(y')=x' }}\lambda_{y\to y'}(x),
\end{equation}
has a stationary distribution that is proportional to $\nu$. 
\end{corollary}
\begin{proof}
 The proof is concluded by simply noting that, due to Proposition~\ref{prop:copy_bal}, for every state $x$ and any $1\leq i\leq h$ we have
 \begin{equation*}
  \sum_{x'\in\ZZ^n_{\geq0}}\nu(x)\sum_{\substack{y\to y'\in\R:\\f_i(y)=x, f_i(y')=x' }}\lambda_{y\to y'}(x)=\sum_{x'\in\ZZ^n_{\geq0}}\nu(x')\sum_{\substack{y\to y'\in\R:\\f_i(y)=x', f_i(y')=x }}\lambda_{y\to y'}(x').
 \end{equation*}
 Hence, by summing both sides over $i$ we obtain that $\nu$ is a stationary measure, hence proportional to a stationary distribution by classic theory of Markov chains.
\end{proof}
 Corollary~\ref{cor:copy_bal} implies that the stationary distributions of the full model can be found, up to a normalizing constant, by studying the stationary distributions of the finite continuous time Markov chain with rates \ref{eq:rate}, as long as the copies are node balanced with respect to them. As already mentioned in the introduction, the node balancing of the copies is intimately related to complex balancing, in a way made precise by the results below.

\begin{theorem}\label{thm:any_kinetics}
 Let $(\G,\Lambda)$ be a stochastic reaction system and let $\nu$ be a measure. Then the following are equivalent:
 \begin{enumerate}
  \item\label{item:everyinj} every injective copy of $\G$ is node balanced with respect to $(\Lambda, \nu)$;
  \item\label{item:cb1} $\nu$ is a complex balanced measure;
  \item\label{item:every} every copy of $\G$ is node balanced with respect to $(\Lambda, \nu)$.
 \end{enumerate}
\end{theorem}
\begin{proof}
 We will prove that \eqref{item:everyinj} implies \eqref{item:cb1}, which implies \eqref{item:every}, which implies \eqref{item:everyinj}.
 \begin{description}
  \item[\eqref{item:everyinj}$\implies$\eqref{item:cb1}.] Fix $y\in\C$ and $x\in\ZZ^n_{\geq0}$. If $x\ngeq y$ then $\lambda_{y\to y'}(x)=0$ for all reactions $y\to y'\in\R$ by \eqref{eq:stoch_cond}. Moreover, for all reactions $y'\to y\in\R$ we have $x+y'-y\ngeq y'$ which implies $\lambda_{y'\to y}(x+y'-y)=0$ and again $\lambda_{y\to y'}(x)=0$ by \eqref{eq:stoch_cond}. It follows that \eqref{eq:cb_def} holds. If $x\geq y$ then by Lemma~\ref{lem:inj_copy_exists} there exists an injective copy $f$ of $\G$ with $f(x)=y$, which is node balanced with respect to $(\Lambda, \nu)$. By injectivity, \eqref{eq:nb_def} becomes \eqref{eq:cb_def}. In conclusion, \eqref{eq:cb_def} holds for every $y\in\C$ and every $x\in\ZZ^n_{\geq0}$, and \eqref{item:cb1} is proven.
  \item[\eqref{item:cb1}$\implies$\eqref{item:every}.] Let $f$ be a copy of $\G$ and let $x\in f(\C)$. Then, summing both sides of \eqref{eq:cb_def} over the complexes $y$ with $f(y)=x$ we obtain \eqref{eq:nb_def}, so \eqref{item:every} holds.
  \item[\eqref{item:every}$\implies$\eqref{item:everyinj}.] This is trivially true as injective copies of $\G$ are copies of $\G$.
 \end{description}
\end{proof}
With more assumptions on the form of the kinetics, a finite number of node balanced copies with respect to $(\Lambda,\nu)$ are sufficient to imply complex balancing of $\nu$, as detailed below. 
\begin{theorem}\label{thm:gen_copy_bal}
 Let $(\G,K^S_{\kappa,\theta})$ be a stochastic reaction system and let $\nu$ be a positive stationary measure. Then there exist two positive and finite constants $M_1, M_2$ such that the following are equivalent:
 \begin{enumerate}
  \item\label{item:1} all the injective copies of $\G$ that intersect the cube $[0,M_1]^n$ are node balanced with respect to $(K^S_{\kappa,\theta}, \nu)$;
  \item\label{item:2} a closed irreducible set $\Gamma$ containing a state $x$ with $\|x\|\geq M_2$ is such that all copies $f$ of $\G$ with $f(\C)\subseteq\Gamma$ are node balanced with respect to $(K^S_{\kappa,\theta}, \nu)$;
  \item\label{item:3} $\nu$ is a complex balanced measure;
 \end{enumerate}
\end{theorem}
The proof is given in Section~\ref{sec:proofs}. In the following, if $f$ is  a copy of $\G$ and $v \in \ZZ^n$, we denote by $f + v$ the function defined by $(f+v)(y) = f(y) + v$ for all $y \in \C$. 
Note that for $v \in \ZZ^n_{\ge 0}$, $f+v$ is a copy of $\G$ while this may not always hold for general $v \in \ZZ^n$.  

\begin{theorem}\label{thm:poisson_copy_bal}
 Let $(\G,K^S_{\kappa,\theta})$ be a stochastic product form system and let $\nu$ be of the form \eqref{eq:general_sm}. Then the following are equivalent:
 \begin{enumerate}
  \item\label{item:onecopy} there is an active, injective, node balanced copy of $\G$ with respect to $(K^S_{\kappa,\theta}, \nu)$;
  \item\label{item:cb2} $\nu$ is a complex balanced measure.
 \end{enumerate}
\end{theorem}
\begin{proof}
 In what follows we will use that, by substituting the rate functions with \eqref{eq:general} and $\nu$ with \eqref{eq:general_sm} and by simplifying, we get that for all $y'\to y''\in\R$ and for all $\tilde{x}\in\ZZ^n_{\geq0}$ with $\tilde{x}\geq y'$ 
  \begin{equation}\label{eq:intermediate}
   \nu(\tilde{x})\lambda_{y'\to y''}(\tilde{x})=c^{\tilde{x}}\kappa_{y'\to y''}\prod_{i=1}^n\prod_{j=1}^{x_i-y'_i}\frac{1}{\theta_i(j)}.
  \end{equation} 
 \begin{description}
  \item[\eqref{item:onecopy}$\implies$\eqref{item:cb2}.] Fix $y\in\C$ and let $x$ be such that $f(x)=y$. By injectivity and node balancing of $f$ we have
  \begin{equation*}
   \nu(x) \sum_{y'\in\C\,:\,y\to y'\in\R}  \lambda_{y\to y'}(x) =  \sum_{y'\in\C\,:\,y'\to y\in\R} \nu(x+y'-y) \lambda_{y'\to y}(x+y'-y),
  \end{equation*}
  which by \eqref{eq:intermediate} becomes
  \begin{equation}\label{eq:intermediate2}
   \sum_{y'\in\C\,:\,y\to y'\in\R}  \kappa_{y\to y'} =  \sum_{y'\in\C\,:\,y'\to y\in\R} c^{y'-y} \kappa_{y'\to y},
  \end{equation}
  which is in turn equivalent to \eqref{eq:cb_def} for mass-action kinetics and hence proves that $(\G, K^D_{\kappa})$ is complex balanced. Moreover, by using \eqref{eq:intermediate} again we have that for all $\tilde{x}\in\ZZ^n_{\geq0}$ with $\tilde{x}\geq y$
  \begin{align*}
   \sum_{y'\in\C\,:\,y'\to y\in\R}  \nu(\tilde{x}+y'-y)\lambda_{y'\to y}(\tilde{x}+y'-y) &= \sum_{y'\in\C\,:\,y'\to y\in\R}c^{\tilde{x}+y'-y}\kappa_{y'\to y}\prod_{i=1}^n \prod_{j=1}^{\tilde{x}_i-y_i}\frac{1}{\theta_i(j)}\\
   &=\sum_{y'\in\C\,:\,y\to y'\in\R}c^{\tilde{x}}\kappa_{y\to y'}\prod_{i=1}^n \prod_{j=1}^{\tilde{x}_i-y_i}\frac{1}{\theta_i(j)}\\
   &=\sum_{y'\in\C\,:\,y'\to y\in\R}\nu(\tilde{x})\lambda_{y\to y'}(\tilde{x}),
  \end{align*}
  where the second equality follows from \eqref{eq:intermediate}. Thus, \eqref{item:cb2} is proven.
  \item[\eqref{item:cb2}$\implies$\eqref{item:onecopy}.] The existence of an injective, active copy follows from Lemma~\ref{lem:inj_copy_exists}, and the fact that is node balanced follows from Theorem~\ref{thm:any_kinetics}.
  \end{description}
\end{proof}

Our last result only concerns mass-action models.

\begin{theorem}\label{thm:poisson_copy_bal2}
 Let $(\G,K^S_\kappa)$ be a stochastic mass-action system and let $\nu$ be of the form \eqref{eq:general_sm}. Then the following are equivalent:
 \begin{enumerate}
  \item\label{item:notinjective} there is a copy $f$ of $\G$ such that $f + v$ is node balanced with respect to $(K^S_{\kappa,\theta},\nu)$ for all $v \in \ZZ^n_{\ge 0}$;
  \item\label{item:notinjective_finite} there is an active copy $f$ of $\G$ such that $f + v$ is node balanced with respect to $(K^S_{\kappa,\theta},\nu)$ for all $v \in \Xi$ where $\Xi\subseteq \ZZ^n_{\geq0}$ is such that the only polynomial of degree at most $d=\max \{ \norm{y}_1 : y \to y' \in \R \}$ that vanishes on $\Xi$ is the zero polynomial;
  \item\label{item:cb3} $\nu$ is a complex balanced measure.
 \end{enumerate}
\end{theorem}
\begin{proof}
 We will prove that \eqref{item:notinjective} implies \eqref{item:notinjective_finite}, which implies \eqref{item:cb3}, which implies \eqref{item:notinjective}. Similarly to what done in the proof of Theorem~\ref{thm:poisson_copy_bal}, we first note that by substituting the rate functions with \eqref{eq:sma} and $\nu$ with \eqref{eq:general_sm} and by simplifying, we get that for all $y'\to y''\in\R$ and for all $\tilde{x}\in\ZZ^n_{\geq0}$ with $\tilde{x}\geq y'$ 
  \begin{equation}\label{eq:intermediate3}
   \nu(\tilde{x})\lambda_{y'\to y''}(\tilde{x})=c^{\tilde{x}}\kappa_{y'\to y''}\frac{1}{(x-y')!}.
  \end{equation} 
 \begin{description}
  \item[\eqref{item:notinjective}$\implies$\eqref{item:notinjective_finite}.] This simply follows from the existence of a set $\Xi$ as described. A finite set $\Xi$ can always be constructed, but to prove \eqref{item:notinjective_finite} it is enough to note that $\Xi$ can be chosen as $\ZZ^n_{\geq0}$, since the only polynomial vanishing on $\ZZ^n_{\geq0}$ is the zero polynomial.
  \item[\eqref{item:notinjective_finite}$\implies$\eqref{item:cb3}.] Let $x$ be in $f(\C)$. By definition of node balancing we have 
  \begin{equation*}
    \nu(x+v) \sum_{y\to y'\in\R:f(y) = x}  \lambda_{y\to y'}(x+v) =  \sum_{y'\to y\in\R:f(y) = x} \nu(x+v+y'-y) \lambda_{y'\to y}(x+v+y'-y)
  \end{equation*}
  for all $v\in\Xi$. By multiplying both sides by $(x+v)!$ and by applying \eqref{eq:intermediate3} we get
  \begin{equation*}
    \sum_{y\in\C\,:\,f(y)=x} \frac{(x+v)!}{(x+v-y)!} \left(\sum_{y'\in\C\,:\,y\to y'\in\R}  \kappa_{y\to y'} - \sum_{y'\in\C\,:\,y'\to y\in\R} c^{y'-y}\kappa_{y'\to y}\right)=0.
  \end{equation*}
  Note that the above is a polynomial of degree at most $d$ in the variable $v$ and it vanishes for all $v\in\Xi$. Hence, it is the zero polynomial and as a consequence $(\G,K^D_\kappa)$ is complex balanced with complex balanced steady state $c$. Then, by using \eqref{eq:intermediate3} again we have that for all $\tilde{x}\in\ZZ^n_{\geq0}$ with $\tilde{x}\geq y$
  \begin{align*}
   \sum_{y'\in\C\,:\,y'\to y\in\R}  \nu(\tilde{x}+y'-y)\lambda_{y'\to y}(\tilde{x}+y'-y) &= \sum_{y'\in\C\,:\,y'\to y\in\R}c^{\tilde{x}+y'-y}\kappa_{y'\to y}\frac{1}{(\tilde{x}-y)!}\\
   &=\sum_{y'\in\C\,:\,y\to y'\in\R}c^{\tilde{x}}\kappa_{y\to y'}\frac{1}{(\tilde{x}-y)!}\\
   &=\sum_{y'\in\C\,:\,y'\to y\in\R}\nu(\tilde{x})\lambda_{y\to y'}(\tilde{x}),
  \end{align*}
  where the second equality follows from \eqref{eq:intermediate3}. Thus, \eqref{item:cb2} is proven.  
  \item[\eqref{item:cb3}$\implies$\eqref{item:notinjective}.] This follows from Theorem~\ref{thm:any_kinetics}.
 \end{description}
\end{proof}
We give here an example showing that condition~\ref{item:notinjective} alone does not imply $\nu$ is complex balanced, unless $\nu$ is assumed to be of the form \eqref{eq:general_sm}.
\begin{example}
 Consider the stochastic mass-action system
  \begin{equation*}
  0 \ce{->[\kappa_1]} A ~,~  3A \cee{->[\kappa_2]} 2A, 
  \end{equation*}
  where the reaction rates have been written above the corresponding reactions. Let $\pi$ be the unique stationary distribution, whose support is in $\{m\in\ZZ\,:\,m\geq2\}$. Since the model is a birth and death chain, $\pi$ is detailed balanced: for all $m\geq2$ we have $\pi(m)q(m,m+1)=\pi(m+1)q(m+1,m)$ \cite{norris:markov}. Specifically, we have
  \begin{equation}\label{eq:example}
   \pi(m)\lambda_{0\to A}(m)=\pi(m+1)\lambda_{3A\to 2A}(m+1)
  \end{equation}
  The deterministically modeled system cannot be complex balanced by Theorem~\ref{thm:classic} because it is not weakly reversible. By \cite[Corollary 19]{cappelletti:complex_balanced}, $\pi$ cannot be complex balanced. However, it follows from \eqref{eq:example} that the copy defined by $f(0)=f(2A)=2$ and $f(A)=f(3A)=3$ is node balanced with respect to $(K^S_\kappa, \pi)$, and the same holds for all the copies $f+v$ with $v\in\ZZ_{\geq0}$. 
\end{example}

\section{Proofs}\label{sec:proofs}

\subsection{Deficiency theory}

The main idea behind the proofs in this section is based upon classical notions of deficiency theory which we will briefly introduce here. We start with giving the definition of deficiency, first given in \cite{feinberg1972}.

\begin{definition}
 The \emph{deficiency} of a reaction network $\G=(\SS, \C, \R)$ is the number $\delta=m-\ell-\dim S$, where $\ell$ is the number of connected components of the directed graph $(\C,\R)$ and $S$ is the stoichiometric subspace of $\G$.
\end{definition}

The following geometric interpretation will be used in our proofs: let $\{e_y\}_{y\in\C}$ be the canonical basis of $\RR^m$, where each coordinate is univocally associated with a complex. Further, define $d_{y\to y'}=e_{y'}-e_{y}$
 for $y\to y'\in\R$. Let $D=\SSann(d_{y\to y'}\,:\, y\to y'\in\R)$. It is proven in \cite{feinberg1972, horn1972general} that $\dim D=m-\ell$. Then, the space $D$ is linearly isomorphic to $S$ if and only if $\delta=0$. Specifically, consider the homomorphism
 \begin{equation}\label{eq:def_phi}
  \begin{array}{rrcl}
    \varphi\colon & \RR^{m} & \to     & \RR^n \\
               & e_y   & \mapsto & y.
  \end{array}
 \end{equation}
 For  $y\to y'\in\R$, we have $\varphi(d_{y\to y'})=y'-y$ and $\varphi_{|D}\colon D\to S$  is thus a surjective homomorphism. Therefore,
 \begin{equation}\label{eq:deficiency_kernel}
  \dim\kernel \varphi_{|D}=\dim D-s=m-\ell-s=\delta,
 \end{equation}
 which implies that $\varphi_{|D}$ is an isomorphism if and only if $\delta=0$. It further follows that the deficiency is a non-negative number.

 The following classical result is proven in \cite{feinberg1972, horn1972general}.
\begin{theorem}\label{thm:classical}
 Let $(\G, K_{\kappa}^D)$ be a deterministic mass-action system. If $\G$ is weakly reversible and its deficiency is zero, then $(\G, K_{\kappa}^D)$ is complex balanced for any choice of rate constants.
\end{theorem}

\subsection{A new model}

Consider a stochastic reaction system $(\G,K_{\kappa,\theta}^S)$ with $\G=(\SS, \C, \R)$ and $K_{\kappa,\theta}^S$ as in \eqref{eq:general}. We will follow an idea first proposed in \cite{cappelletti:complex_balanced}, and consider the reaction network $\tilde{G}=(\tilde{\SS}, \tilde{\C}, \tilde{\R})$ with
\begin{equation}\label{eq:gtilde_def}
 \tilde{\SS}=\SS\cup\{A_y\,:\,y\in\C\},\quad \tilde{\C}=\{y+A_y\,:\,y\in\C\},\quad\text{and}\quad \tilde{\R}=\{y+A_y\to y'+A_{y'}\,:\,y\to y'\in\R\}.
\end{equation}
It follows from the construction above that $\tilde{\G}$ is weakly reversible if and only if $\G$ is weakly reversible. Moreover, the following holds.
\begin{lemma}\label{lem:def0}
 The deficiency of $\tilde{\G}$ as defined in \eqref{eq:gtilde_def} is zero.
\end{lemma}
\begin{proof}
 Consider the linear homomorphism $\varphi$ as defined in \eqref{eq:def_phi}, for the reaction network $\tilde{G}$. Since the vectors associated with the complexes $y+A_y$ are linear independent, $\varphi$ is an isomorphism and as a consequence the deficiency of $\tilde{G}$ is 0.
\end{proof}
We associate $\tilde{G}$ with a general mass-action kinetics $K^S_{\kappa, \tilde{\theta}}$ such that the rate of any reaction $y+A_y\to y'+A_{y'}\in\tilde{\R}$ is of the form
\begin{equation*}
K^S_{\kappa, \tilde{\theta}}(y+A_y\to y'+A_{y'})(x,u)=K^S_{\kappa, \theta}(y\to y')(x)u_{A_y}\quad\text{for all }(x,u)\in\ZZ^n_{\geq0}\times\ZZ^m_{\geq0}.
\end{equation*}
The following is the key result of this paper, on which the proof of the other main theorems are based on.
\begin{theorem}\label{thm:key}
 Let $(\G,K_{\kappa,\theta}^S)$ be a stochastic reaction system with $K_{\kappa,\theta}^S$ as in \eqref{eq:general}. Let $\nu$ be a $\sigma$-finite measure on $\ZZ^n_{\geq0}$. Define $\Upsilon$ as the set of pairs $(x,y)\in\ZZ^n_{\geq0}\times\C$ such that there exists an active, injective copy $f$ of $\G$ that is node balanced with respect to $(K_{\kappa,\theta}^S,\nu)$, satisfies $f(y)=x$, and fulfils $f(\C)\subseteq\supp\nu$. Let $\hat{\R}$ be a subset of reactions whose reaction vectors form a basis of $S$. Assume that
 \begin{enumerate}
  \item there exists $\hat{x}\in\ZZ^n_{\geq0}$ such that $(\hat{x},y)\in\Upsilon$ for all $y\in\C$;
  \item for all $y^\star\to y^{\star\star}\in\R\setminus\hat{\R}$, there exist a sequence of (potentially repeated) reactions $\{y_i\to y_i'\}_{i=1}^h$ contained in $\{y^\star\to y^{\star\star}\}\cup\hat{\R}$ such that
  \begin{enumerate}
   \item $y^\star\to y^{\star\star}=y_1\to y_1'$;
   \item $\sum_{i=1}^h (y_i'-y_i)=0$;
   \item for all $j=0,1,2,\dots,h-1$ we have
   \begin{equation*}
    \left(\hat{x}+\sum_{i=1}^j (y_i'-y_i), y_{j+1}\right)\in\Upsilon.
   \end{equation*}
  \end{enumerate}
 \end{enumerate}
 Then, the deterministic mass-action system $(\G,K^D_\kappa)$ is complex balanced.
\end{theorem}
\begin{proof}
 Consider $(x,e_y)\in\ZZ^n_{\geq0}\times\ZZ^m_{\geq0}$. If $x\geq y$, then any reaction of the form $y+A_y\to y'+A_{y'}$ can take place and move the process from $(x,e_y)$ to $(x+y'-y, e_{y'})$. From $x\geq y$ it follows that $x+y'-y\geq y'$ hence at $(x+y'-y, e_{y'})$ any reaction of the form $y'+A_{y'}\to y''+A_{y''}$ can take place and move the process to $(x+y''-y, e_{y''})$, and so on. Since $\G$ is weakly reversible it is possible to eventually return to the original state $(x,e_y)$, and the closed irreducible set containing $(x,e_y)$ is precisely given by $\{(x+\tilde{y}-y,e_{\tilde{y}})\,:\,y\in\C_y\}$, where $\C_y$ denotes the connected component of $(\C,\R)$ containing $y$. We denote the closed irreducible set containing $(x,e_y)$ by $\Gamma_{(x,e_y)}$. Of course the notation is not bijective: for example if $x\geq y$ then $\Gamma_{(x,e_y)}=\Gamma_{(x+y'-y,e_{y'})}$. If $(x,y)\in\Upsilon$, then be definition there exists an active, injective copy $f$ of $\G$ that is node balanced with respect to $(K_{\kappa,\theta}^S,\nu)$ and that satisfies $f(y)=x$. Hence, by definition of node balancing, for all $\tilde{y}\in\C_y$ we have
 \begin{equation*}
  \sum_{y'\in\C: \tilde{y}\to y'\in\R} \nu(x+\tilde{y}-y)\lambda_{\tilde{y}\to y'}(x+\tilde{y}-y)=\sum_{y'\in\C: y'\to\tilde{y}\in\R} \nu(x+y'-y)\lambda_{y'\to\tilde{y}}(x+y'-y),
 \end{equation*}
 which implies
 \begin{multline*}
  \sum_{y'\in\C: \tilde{y}+A_{\tilde{y}}\to y'+A_{y'}\in\tilde{\R}} \nu(x+\tilde{y}-y)\lambda_{\tilde{y}+A_{\tilde{y}}\to y'+A_{y'}}(x+\tilde{y}-y,e_{\tilde{y}})\\
  =\sum_{y'\in\C: y'+A_{y'}\to\tilde{y}+A_{\tilde{y}}\in\tilde{\R}} \nu(x+y'-y)\lambda_{y'+A_{y'}\to\tilde{y}+A_{\tilde{y}}}(x+y'-y,e_{y'}).
 \end{multline*}
 Since $\Gamma_{(x,e_y)}$ is finite, it follows from the equation above that $\nu$ restricted to $\Gamma_{(x,e_y)}$ is proportional to the unique stationary distribution of $(\tilde{\G}, K^S_{\kappa, \tilde{\theta}})$ with support $\Gamma_{(x,e_y)}$. Since $\tilde{\G}$ has deficiency zero by Lemma~\ref{lem:def0} and is weakly reversible because $\G$ is weakly reversible, by Theorems~\ref{thm:classical} and \ref{thm:classic} we have
 \begin{equation}\label{eq:key}
  \nu(x+\tilde{y}-y)=M_{(x,e_y)}\tilde{c}_{A_y}c^x\prod_{j=1}^{x_i}\frac{1}{\theta_i(j)},
 \end{equation}
 where $(c,\tilde{c})\in\RR^n_{>0}\times\RR^m_{>0}$ is a complex balanced steady state of $(\tilde{\G}, K^D_\kappa)$ and $M_{(x,e_y)}$ is a proportionality constant depending on the closed irreducible set $\Gamma_{(x,e_y)}$, hence the notation is not bijective and $M_{(x,e_y)}=M_{(x+\tilde{y}-y,e_{\tilde{y}})}$ for all $\tilde{y}\in\C_y$.
 
 We now want to show that there exists a vector $\hat{c}\in\RR_{>0}^n$ such that
 \begin{equation}\label{eq:toprove}
  \frac{\tilde{c}_{A_y'}}{\tilde{c}_{A_y}}=\hat{c}^{y'-y}\quad\text{for all }y\to y'\in \R.
 \end{equation}
 First, note that \eqref{eq:toprove} holds for all $y\to y'\in \hat{\R}$, because it is equivalent to
 $$(y'-y)^\top \log\hat{c}=\log\tilde{c}_{A_y'}-\log\tilde{c}_{A_y}\quad\text{for all }y\to y'\in \hat{\R}.$$
 The latter has a solution because the reaction vectors of the reactions in $\hat{R}$ are linearly independent. If $y\to y'\in\R\setminus\hat{\R}$ then by hypothesis there exists a sequence of reactions $\{y_i\to y_i'\}_{i=1}^h$ contained in $\{y\to y'\}\cup\hat{\R}$ such that $y_1\to y_1'=y\to y'$ and 
 \begin{equation*}
  \left(\hat{x}+\sum_{i=1}^j (y_i'-y_i), y_{j+1}\right)\in\Upsilon.
 \end{equation*}
 for all $j=0,1,2,\dots, h-1$. Hence, it follows from $y'-y+\sum_{i=1}^h (y_i'-y_i)=0$ and from applying \eqref{eq:key} recursively that 
 \begin{align*}
  M_{(\hat{x},y_1)}&=M_{(\hat{x}+y'-y+\sum_{i=1}^h (y_i'-y_i),y_1)}\\
  &=\frac{\tilde{c}_{A_{y_{h-1}'}}}{\tilde{c}_{A_{y_1}}}M_{(\hat{x}+y'-y+\sum_{i=1}^h (y_i'-y_i),y_{h-1}')}=\frac{\tilde{c}_{A_{y_{h-1}'}}}{\tilde{c}_{A_{y_1}}}M_{(\hat{x}+y'-y+\sum_{i=1}^{h-1} (y_i'-y_i),y_{h-1})}\\ 
  &=\frac{\tilde{c}_{A_{y_{h-1}'}}\tilde{c}_{A_{y_{h-2}'}}}{\tilde{c}_{A_{y_1}}\tilde{c}_{A_{y_{h-1}}}}M_{(\hat{x}+y'-y+\sum_{i=1}^{h-1} (y_i'-y_i),y_{h-2}')}=\frac{\tilde{c}_{A_{y_{h-1}'}}\tilde{c}_{A_{y_{h-2}'}}}{\tilde{c}_{A_{y_1}}\tilde{c}_{A_{y_{h-1}}}}M_{(\hat{x}+y'-y+\sum_{i=1}^{h-2} (y_i'-y_i),y_{h-2})}\\
  &=\dots=M_{(\hat{x},y_1)}\prod_{i=1}^h \frac{\tilde{c}_{A_{y_i'}}}{\tilde{c}_{A_{y_i}}}.
 \end{align*}
 As a consequence,
 \begin{equation*}
  \prod_{i=1}^h \frac{\tilde{c}_{A_{y_i'}}}{\tilde{c}_{A_{y_i}}}=1.
 \end{equation*}
 If $y\to y'$ appears $\alpha$ times in the sequence $\{y_i\to y_i'\}_{i=1}^h$, and all other reactions $\tilde{y}\to \tilde{y}'$ appear $\beta_{\tilde{y}\to \tilde{y}'}$ times, then
 \begin{equation*}
  \frac{\tilde{c}_{A_{y'}}}{\tilde{c}_{A_{y}}}=\prod_{\tilde{y}\to \tilde{y}'\in\hat{\R}} \left(\frac{\tilde{c}_{A_{\tilde{y}}}}{\tilde{c}_{A_{\tilde{y}'}}}\right)^{\frac{\beta_{\tilde{y}\to \tilde{y}'}}{\alpha}}=\hat{c}^{-\frac{1}{\alpha}\sum_{\tilde{y}\to \tilde{y}'\in\hat{\R}}\beta_{\tilde{y}\to \tilde{y}'}(y'-y)}=\hat{c}^{y'-y}.
 \end{equation*}
 Hence, \eqref{eq:toprove} is proven. For all $y\in\C$ we have $(\hat{x},y)\in\Upsilon$, which implies
 \begin{equation*}
  \sum_{y'\in\C: y\to y'\in\R} \nu(\hat{x})\lambda_{y\to y'}(\hat{x})=\sum_{y'\in\C: y'\to y\in\R} \nu(x+y'-y)\lambda_{y'\to\tilde{y}}(x+y'-y).
 \end{equation*}
 By substituting the expressions for the transition rates and $\nu$ (over the image of node balanced copies) and by simplifying we obtain
 \begin{equation*}
  \sum_{y'\in\C: y\to y'\in\R} c^x\kappa_{y\to y'}=\sum_{y'\in\C: y'\to y\in\R} \frac{\tilde{c}_{A_{y'}}}{\tilde{c}_{A_y}}c^{x+y'-y}\kappa_{y'\to y}
 \end{equation*}
 which by \eqref{eq:toprove} implies
  \begin{equation*}
  \sum_{y'\in\C: y\to y'\in\R} \kappa_{y\to y'}\prod_{i=1}^n (c_i \hat{c}_i)^{y_i}=\sum_{y'\in\C: y'\to y\in\R} \kappa_{y'\to y}\prod_{i=1}^n (c_i \hat{c}_i)^{y_i'},
 \end{equation*}
 hence a positive complex balanced steady state exists for $(\G,K^D_\kappa)$, which concludes the proof.
\end{proof}

\subsection{Proof of Theorem~\ref{thm:cbiffcb}}

By Theorem~\ref{thm:productform} we already know that if $(\G,K^D_\kappa)$ is complex balanced, then there exists a $\sigma$-finite, positive, stationary measure for $(\G,K^S_{\kappa, \theta})$ given by $\nu$ as in \eqref{eq:general_sm}. Checking that $\nu$ is complex balanced is not difficult, since by substituting the reaction rates with \eqref{eq:general} and $\nu$ with \eqref{eq:general_sm}, equation \eqref{eq:cb_def} simplifies to
\begin{equation*}
 \sum_{y'\in\C\,:\,y\to y'\in\R}  \kappa_{y\to y'} =  \sum_{y'\in\C\,:\,y'\to y\in\R} c^{y'-y} \kappa_{y'\to y},
\end{equation*}
which holds because $(\G,K^D_\kappa)$ is complex balanced.

Conversely, assume a $\sigma$-finite, positive, complex balanced measure for $(\G,K^S_{\kappa, \theta})$ exists. Then, all copies of $\G$ are node balanced with respect to $(\G,K^S_{\kappa, \theta})$ by Theorem~\ref{thm:any_kinetics}. Moreover, since $\nu$ is positive then for all copies $f$ of $\G$ we have $f(\C)\subseteq\supp \nu$. Hence, the proof is concluded by Theorem~\ref{thm:key} by choosing $\hat{x}\in\ZZ^n_{\geq0}$ large enough.

\subsection{Proof of Theorem~\ref{thm:gen_copy_bal}}

It follows from Theorem~\ref{thm:any_kinetics} that \eqref{item:3} implies \eqref{item:1} for any $M_1$. Moreover, the existence of a positive complex balanced measure implies that $\G$ is weakly reversible \cite[Theorem 4.4]{CJ2018} hence $\ZZ^n_{\geq0}$ is union of closed irreducible sets\cite{craciun:dynamical}. As a consequence, the existence of a closed irreducible set as in \eqref{item:2} is implied by Theorem~\ref{thm:any_kinetics} for any $M_2$.

Conversely, since $\nu$ is positive, it is always possible to choose $M_1$ large enough such that all injective copies intersecting $[0,M_1]^n$ being node balanced with respect to $(K^S_{\kappa,\theta})$ implies that the assumptions of Theorem~\ref{thm:key} are satisfied. Hence, for large enough $M_1$ \eqref{item:1} implies that $(\G,K^D_\kappa)$ is complex balanced. In turn, this implies that $\nu$ is proportional to \eqref{eq:general_sm} on every closed irreducible set because the model is non-explosive by Theorem~\ref{thm:finite_non_explosive} (hence all stationary measures are proportional to each other on every closed irreducible set by standard Markov chain theory \cite{norris:markov}). By substituting the reaction rates with \eqref{eq:general} and $\nu$ as proportional to \eqref{eq:general_sm}, equation \eqref{eq:cb_def} simplifies to
\begin{equation*}
 \sum_{y'\in\C\,:\,y\to y'\in\R}  \kappa_{y\to y'} =  \sum_{y'\in\C\,:\,y'\to y\in\R} c^{y'-y} \kappa_{y'\to y},
\end{equation*}
which holds because $(\G,K^D_\kappa)$ is complex balanced.

Similarly, since $\nu$ is positive, all injective copies in a closed irreducible set with large enough states being node balanced with respect to $(K^S_{\kappa,\theta})$ implies that the assumptions of Theorem~\ref{thm:key} are satisfied. By following the same reasoning as above, this implies \eqref{item:3} and concludes the proof.

\section*{Acknowledgements}

DC was supported by the MIUR grant ‘Dipartimenti di Eccellenza 2018-2022’ (E11G18000350001).
 
\bibliographystyle{plain}
\bibliography{bib}

\end{document}